\documentstyle[12pt]{article}
\newcommand{\nc}{\newcommand}  \nc{\ov}{\over} \nc{\cd}{\cdots}
\nc{\ra}{\rightarrow} \nc{\iy}{\infty}
\nc{\be}{\begin{equation}} \nc{\ee}{\end{equation}}
\nc{\bZ}{{\bf Z}} \nc{\dl}{\delta} \nc{\la}{\lambda}
\nc{\inv}{^{-1}}  \nc{\ph}{\phi} \nc{\ps}{\psi}
\nc{\twotwo}[4]{\left(\begin{array}{cc}#1&#2\\#3&#4\end{array}\right)}
\nc{\twoone}[2]{\left(\begin{array}{c}#1\\#2\end{array}\right)}
\nc{\cP}{{\cal {P}}} \renewcommand{\sp}{\vskip1ex} \nc{\noi}{\noindent}
\begin{document} 

\begin{center}{\bf ON A TOEPLITZ DETERMINANT IDENTITY\\
OF BORODIN AND OKOUNKOV}\end{center}\vskip.1ex
\begin{center}{\bf Estelle L. Basor\footnote{Supported by National Science 
Foundation grant DMS--9970879.} and Harold Widom\footnote{Supported by National 
Science Foundation grant DMS--9732687.}}\end{center}\vskip.1ex

\begin{quotation} {\small In this note we give two other proofs of an
identity of A.~Borodin and A.~Okounkov
which expresses a Toeplitz determinant in terms of the Fredholm determinant of a
product of two Hankel operators. The second of these proofs yields a generalization
of the identity to the case of block Toeplitz determinants.}\end{quotation}\vskip.2ex

The authors of the title proved in [2] an elegant identity expressing a Toeplitz
determinant in terms of the Fredholm determinant of an infinite matrix which (although not
described as such) is the product of two Hankel matrices. The proof used
combinatorial theory, in particular a theorem of Gessel expressing a Toeplitz
determinant as a sum over partitions of products of Schur functions. The purpose
of this note is to give two other proofs of the identity. The first uses an identity
of the second author [4] for the quotient of Toeplitz determinants in which the same
product of Hankel matrices appears and the second, which is more direct and 
extends the identity to the case of block Toeplitz determinants, consists
of carrying the first author's collaborative proof [1] of the strong Szeg\"o limit 
theorem one  step further.

We begin with the statement of the identity of [2], changing notation slightly.
If $\ph$ is a function on the unit circle with Fourier coefficients $\ph_k$
then $T_n(\ph)$ denotes the Toeplitz matrix $(\ph_{i-j})_{i,j=0,\cd,n-1}$ and $D_n(\ph)$
its determinant. Under general conditions $\ph$ has a representation
$\ph=\ph_+\,\ph_-$ where $\ph_+$ (resp. $\ph_-$) extends to a nonzero analytic
function in the interior (resp. exterior) of the circle. We assume
that $\ph$ has geometric mean 1, and normalize $\ph_{\pm}$ so that
$\ph_+(0)=\ph_-(\iy)=1$. Define the infinite matrices $U_n$ and $V_n$ acting
on $\ell^2(\bZ^+)$, where $\bZ^+=\{0,\,1,\cd\}$, by
\[U_n(i,j)=(\ph_-/\ph_+)_{n+i+j+1},\ \ \ V_n(i,j)=(\ph_+/\ph_-)_{-n-i-j-1}\]
and the matrix $K_n$ acting on $\ell^2(\{n,\,n+1,\cd\})$ by
\[K_n(i,j)=\sum_{k=1}^{\iy}(\ph_-/\ph_+)_{i+k}\;(\ph_+/\ph_-)_{-k-j}.\]
Notice that $K_n$ becomes $U_nV_n$ under the obvious identification of
$\ell^2(\{n,\,n+1,\cd\})$ with $\ell^2(\bZ^+)$. It is easy to check that, aside from
a factor $(-1)^{i+j}$ which does not affect its Fredholm determinant, the entries
of $K_n$ are the same as given by the integral formula (2.2) of [2]. The
formula of Borodin and Okounkov is
\be D_n(\ph)=Z\,\det\,(I-K_n),\label{BOform}\ee
where
\[Z=\exp\,\{\sum_{k=1}^{\iy}k\,(\log\,\ph)_k\,(\log\,\ph)_{-k}\}=\lim_{n\ra\iy}D_n(\ph).\]
(The last identity is the strong Szeg\"o limit theorem.) This identity is especially
useful for obtaining refined asymptotic results as $n\ra\iy$.

Two versions of (\ref{BOform}) were proved in [2]. One was algebraic
and was an identity of formal power series and the other was analytic and assumed that
the regions of analyticity of $\ph_{\pm}$ included neighborhoods of the unit circle
although, as the authors point out, an approximation argument can be used to extend the
range of validity. The requirements for our proofs are that
$\log\,\ph_{\pm}$ be bounded and
$\sum_{k=-\iy}^{\iy}|k|\,|(\log\,\ph)_k|^2<\iy$. \footnote{The bounded functions
$f$ satisfying $\sum_{k=-\iy}^{\iy}|k||f_k|^2<\iy$
form a Banach algebra under a natural norm and for any such $f$ the Hankel
matrix $(f_{i+j})$ acting on $\ell^2(\bZ^+)$ is Hilbert-Schmidt. Thus if $\log\ph_{\pm}$
belong to this algebra so do $\ph_-/\ph_+$ and $\ph_+/\ph_-$ and it follows that $U_n$ and $V_n$
are Hilbert-Schmidt so $K_n$ is trace class. Moreover the Szeg\"o limit
theorem holds for such $\ph$. See [5] or, for this and a lot more, [3].}\sp

\begin{center}{\bf First proof}\end{center}

To state the relevant result of [4] we define the vectors $U_n\dl$
and $V_n\dl$ in $\bZ^+$ by
\[U_n\dl(i)=(\ph_-/\ph_+)_{n+i},\ \ \ V_n\dl(i)=(\ph_+/\ph_-)_{-n-i}.\]
(These are not the results of acting on a vector $\dl$ by the operators $U_n$ and $V_n$
since $-1\not\in\bZ^+$, but the notation suggests this.) The result is the following proposition. \sp

If $I-U_n\,V_n$ is invertible then so is $T_n(\ph)$ and
\be{D_{n-1}(\ph)\ov D_n(\ph)}=1-\left((I-U_n\,V_n)\inv U_n\dl,\,V_n\dl\right),\label{Wform}\ee
where the inner product denotes the sum of the products of the components.\sp

The formula appears on p. 341 of [4] in different notation.
To derive (\ref{BOform}) from this we assume temporarily that $I-U_n\,V_n$ is
invertible for all $n$ (and therefore so is $I-V_n\,U_n$) and compute the upper-left
entry of $(I-V_{n-1}\,U_{n-1})\inv$
in two different ways. This entry equals the upper-left entry of $(I-K_{n-1})\inv$,
and Cramer's rule says that the inverse of the entry equals
\[{\det\,(I-K_{n-1})\ov\det\,(I-K_n)}.\]
On the other hand, there is a general formula which says that if one has a
$2\times 2$ block matrix
$\twotwo{A}{B}{C}{D}$ then the upper-left block of its inverse equals $(A-BD\inv C)\inv$.
Here $A$ and $D$ are square and the various inverses are assumed to exist.
In our case the
large matrix is $I-V_{n-1}\,U_{n-1}$ and $A$ is $1\times1$. It is easy to
see that
\[A=1-(V_n\dl,\,U_n\dl),\ \ \ D=I-V_n\,U_n,\ \ \ C=-V_n\,U_n\dl,\ \ \ B=-U_n\,V_n\dl,\]
the last interpreted as a row vector. The formula says that the inverse of
the upper-left entry of the inverse equals
\[1-(V_n\dl,\,U_n\dl)-\left((I-V_n U_n)\inv \,V_n\,U_n\dl,\,U_n\,V_n\dl\right)\]
\[=1-(V_n\dl,\,U_n\dl)-\left(U_n\,(I-V_n U_n)\inv\, V_n\,U_n\dl,\,V_n\dl\right)\]
\[=1-(V_n\dl,\,U_n\dl)-\left(\left[(I-U_n V_n)\inv-I\right]\,U_n\dl,\,V_n\dl\right)\]
\[=1-\left((I-U_n V_n)\inv U_n\dl,\,V_n\dl\right),\]
which is the right side of (\ref{Wform}). Thus we have established
\[{D_{n-1}(\ph)\ov D_n(\ph)}={\det\,(I-K_{n-1})\ov\det\,(I-K_n)},\]
which shows that (\ref{BOform}) holds for some constant $Z$. That $Z$ is
as stated follows by letting $n\ra\iy$.

To remove the restriction that $I-U_n\,V_n$ be invertible for all $n$, we
introduce a complex parameter $\la$ and replace $\ph$ by $\ph^\la=\exp\,(\la\log\ph)$.
Then both sides of (\ref{BOform}) are entire functions of $\la$ and are equal
when $\la$ is so small that $\parallel \!\ph_-^\la/\ph_+^\la-1\!\parallel_{\iy}<1$
and $\parallel \!\ph_+^\la/\ph_-^\la-1\!\parallel_{\iy}<1$, for then all $U_n$ and $V_n$
have operator norm less than 1 so all $I-U_n\,V_n$ are invertible.
Since the two sides sides of (\ref{BOform}) are equal for small $\la$ they are equal
for all $\la$.\sp

\begin{center}{\bf Second proof}\end{center}

Denote by $T(\phi)$ the semi-infinite Toeplitz matrix
$(\ph_{i-j})_{i,\,j\ge0}$. Then $T(\phi_{-})$ and $T(\phi_{+})$ are upper-triangular 
and lower-triangular resepectively. It follows that if $P_{n}$ is the diagonal matrix 
whose first
$n$ diagonal entries are all 1 and whose other entries are 0 then
\[P_{n}\,T(\phi_{+}) = P_{n}\,T(\phi_{+})\,P_{n},\ \ \ 
T(\phi_{-})\,P_{n} = P_{n}\,T(\phi_{-})\,P_{n}.\]
Observe that $T_n(\ph)$ is the upper-left $n\times n$ block of $P_n\,T_n(\ph)\,P_n$.
Using the above, we can write \footnote{It is an easy general fact that if 
$\psi_1\in \overline{H^{\iy}}$ or $\psi_2\in H^{\iy}$ than $T(\psi_1\psi_2)
=T(\psi_1)T(\psi_2)$. In particular $T(\phi_{\pm})$ are invertible with inverses
$T(\phi_{\pm}^{-1})$. Recall that $H^{\iy}$ consists of all $\psi\in L^{\iy}$ such that
$\psi_k=0$ when $k<0$.}
\[P_n\,T(\ph)\,P_n
=P_{n}\,T(\phi_{+})\,T(\phi_{+}^{-1})\,T(\phi)\,T(\phi_{-}^{-1})\,T(\phi_{-})\,P_{n}\]
\[=P_{n}\,T(\phi_{+})\,P_{n}\,T(\phi_{+}^{-1})\,T(\phi)\,
T(\phi_{-}^{-1})\,P_{n}\,T(\phi_{-})\,P_{n}.\]

Now the upper-left blocks of $P_{n}\,T(\phi_{\pm})\,P_{n}$ are $T_n(\ph_\pm)$, which
are triangular matrices with diagonal entries all 1, by our assumed normalization.
Therefore they have determinant one, so $D_n(\ph)$ equals the determinant of the 
upper-left block of
$P_{n}\,T(\phi_{+}^{-1})\,T(\phi)\,T(\phi_{-}^{-1})\,P_{n}$.
Set 
\[T(\phi_{+}^{-1})\,T(\phi)\,T(\phi_{-}^{-1}) = A.\]
Then the determinant of the upper-left block of $P_nAP_n$
equals $\det\,(P_nAP_n+Q_n)$, where $Q_{n}= I -P_{n}.$ Now $A$ is invertible 
and differs from $I$ by a trace class operator (we shall see this
in a moment). Therefore
\[\det\,( P_{n}A P_{n}+Q_{n})=\det A \,\det\,(A^{-1} P_{n}A P_{n}+A^{-1}Q_{n})\]
 \[= \det A \,\det\,(A^{-1}( I-Q_{n})A P_{n}+A^{-1}Q_{n}) 
 = \det A \,\det\,(P_{n}-A^{-1}Q_{n}A P_{n}+A^{-1}Q_{n})\]
 \[= \det A \,\det\,(P_{n} + A^{-1}Q_{n})\,\det\,(I - Q_{n}A P_{n}),\]
 since  $P_{n}Q_n=0$. The determinant of the operator on the right equals one,
 again since  $P_{n}Q_n=0$. Moreover
 \[\det\,(P_n+A\inv Q_n)=\det\, (I-(I-A\inv)Q_n)=\det \,(I-Q_n(I-A\inv)Q_n).\]
 We have shown 
 \be D_n(\ph)=\det A \, \det \,(I-Q_n(I-A\inv)Q_n).\label{Dnform}\ee
 
 It remains to show that this is the same as (\ref{BOform}). First, $A$ is similar 
 via the invertible operator $T(\phi_+)$ to 
$T(\phi)\,T(\phi_{-}^{-1})\,T(\phi_{+}^{-1}).$ Therefore
\be\det A=\det\,T(\phi)\,T(\phi_{-}^{-1})\,T(\phi_{+}^{-1})=\det\,T(\ph)\,T(\ph\inv).
\label{detA}\ee
This is a representation of the constant $Z$ in the strong Szeg\"o limit theorem
[1, 3, 5].
 
 Next
 \be A\inv=T(\ph_-)\,T(\ph)\inv\,T(\ph_+)=T(\ph_-)\,T(\ph_+\inv)\,T(\ph_-\inv)\,T(\ph_+)
 =T(\ph_-/\ph_+)\,T(\ph_+/\ph_-).\label{Ainv}\ee
 Because $\ph_-/\ph_+$ and $\ph_+/\ph_-$ are reciprocals we see that the $i,j$ entry
 of this matrix equals
\[\delta_{i,j}-\sum_{k=1}^{\iy}(\ph_-/\ph_+)_{i+k}\,(\ph_+/\ph_-)_{-k-j},\]
and so $\det\,(I-Q_n(I-A\inv)Q_n)$ equals $\det\,(I-K_n)$. (This also shows
that $A\inv$ differs from $I$ by a trace class operator, so the same is true of $A$.)
This gives (\ref{BOform}) with $Z=\det T(\ph)\,T(\ph\inv)$.\sp 

Let us see how to modify this argument for the case of block Toeplitz determinants, where
$\ph$ is a matrix-valued function. We assume the factorization
$\ph=\ph_+\ph_-$, the order of the factors being important now, where $\ph_{\pm}^{\pm1}$ 
belong to the algebra described in footnote 3
and $\ph_+^{\pm1}\in H^{\iy},\ \ph_-^{\pm1}\in\overline{H^{\iy}}$. Then (\ref{Dnform})
is derived without change as is formula (\ref{detA}) for $\det A$ since
$\ph\inv=\ph_-\inv\ph_+\inv$. But (\ref{Ainv}) no longer holds because it 
would require $\ph=\ph_-\ph_+$, which does not hold. But if we also assume
a factorization $\ph=\psi_-\psi_+$, with $\psi_{\pm}$ having properties analogous to 
those of $\ph_{\pm}$, we can replace (\ref{Ainv}) by
\[ A\inv=T(\ph_-)\,T(\psi_+\inv)\,T(\psi_-\inv)\,T(\ph_+)
 =T(\ph_-\psi_+\inv)\,T(\psi_-\inv\ph_+).\]
 Now $\ph_-\psi_+\inv$ and $\psi_-\inv\ph_+$ are mutual inverses and we deduce that in
 this case (\ref{BOform}) holds with $Z=\det T(\ph)\,T(\ph\inv)$ 
 and $K_n$ the matrix,
 thought of as acting on $\ell^2(\{n,\,n+1,\cd\})$, with $i,j$ entry
 \[\sum_{k=1}^{\iy}(\ph_-\psi_+\inv)_{i+k}\,(\psi_-\inv\ph_+)_{-k-j}.\]

\sp \begin{center}{\bf Acknowledgement}\end{center}
 
The authors thank Alexander Its who, after seeing [2] and the 
first proof presented here, asked the first author whether there was a direct proof. 
The second proof was the result.

\vspace{4ex}

\noi {\large \bf  References} \vspace{2ex}

\noi[1]\ \ E. Basor and J. W. Helton, {\it A new proof of the Szeg\"o limit
theorem and new results for Toeplitz operators with discontinuous symbol},
J. Operator Th. {\bf 3} (1980) 23--39.\vspace{1.5ex}

\noi[2]\ \ A. Borodin and A. Okounkov, {\it A Fredholm determinant formula for
Toeplitz determinants}, preprint, math.CA/9907165.\vspace{1.5ex}

\noi[3]\ \ A. B\"ottcher and B. Silbermann, Analysis of Toeplitz Operators,
Akademie-Verlag, Berlin, 1989.\vspace{1.5ex}

\noi[4]\ \ H. Widom, {\it Toeplitz determinants with singular generating functions},
Amer.~J.~Math. {\bf 95} (1973) 333--383.\vspace{1.5ex}

\noi[5]\ \ H. Widom, {\it Asymptotic behavior of block Toeplitz matrices
and determinants II}, Adv. Math. {\bf 21} (1976) 1--29.\vspace{3ex}

\hspace{-2em}\begin{tabular}{lll} Department of Mathematics&&Department of Mathematics
\\California Polytechnic State University&&University of California\\
San Luis Obispo, CA 93407 USA&&Santa Cruz, CA 95064 USA\\
ebasor@calpoly.edu&&widom@math.ucsc.edu\end{tabular}\sp
\vspace{3ex}


\end{document}